\documentclass[a4paper,12pt]{article}
\usepackage{tabularx} 
\usepackage{amsmath,amsthm,verbatim,amssymb,amsfonts,amsopn,amscd}
\usepackage{graphicx} 
\usepackage[margin=1in,a4paper]{geometry} 
\usepackage[utf8]{inputenc}
\usepackage[T1]{fontenc}

\usepackage[UKenglish]{babel}
\usepackage{enumerate}
\usepackage{cite} 
\usepackage[final]{hyperref} 
\hypersetup{
	colorlinks=true,       
	linkcolor=blue,        
	citecolor=blue,        
	filecolor=magenta,     
	urlcolor=blue         
}

\usepackage{tikz-cd}
\tikzcdset{row sep/normal=2.2em}
\tikzcdset{arrow style=math font}

\theoremstyle{definition}
\newtheorem{definition}{Definition}
\newtheorem{remark}{Remark}

\theoremstyle{plain}
\newtheorem{theorem}{Theorem}
\newtheorem{corollary}{Corollary}
\newtheorem{lemma}{Lemma}
\newtheorem{proposition}{Proposition}

\newcommand{\set}[2]{\{{#1}\,\mid\,{#2}\}}

\DeclareMathOperator{\calA}{\mathcal{A}}

\DeclareMathOperator{\DN}{DN}

\DeclareMathOperator{\D}{D}

\newcommand{\IN}{\mathbb{N}}

\newcommand{\IQ}{\mathbb{Q}}
\newcommand{\IR}{\mathbb{R}}

\newcommand{\into}{\hookrightarrow}

\DeclareMathOperator{\Coo}{C^\infty}

\newcommand{\la}{\langle}
\newcommand{\ra}{\rangle}

\title{Affine connections and second-order affine structures}
\author{\href{mailto:filip.bar@cantab.net}{Filip B\'{a}r}}
\date{\itshape Dedicated to my good friend Tom Rewwer on the occasion of his 35th birthday.}

\begin{document}
	
	\maketitle
	
	\begin{abstract}
		Smooth manifolds have been always understood intuitively as spaces with an affine geometry on the infinitesimal scale. In Synthetic Differential Geometry this can be made precise by showing that a smooth manifold carries a natural structure of an infinitesimally affine space. This structure is comprised of two pieces of data: a sequence of symmetric and reflexive relations defining the tuples of mutual infinitesimally close points, called an infinitesimal structure, and an action of affine combinations on these tuples. For smooth manifolds the only natural infinitesimal structure that has been considered so far is the one generated by the first neighbourhood of the diagonal. In this paper we construct natural infinitesimal structures for higher-order neighbourhoods of the diagonal and show that on any manifold any symmetric affine connection extends to a second-order infinitesimally affine structure.
	\end{abstract}
	
	\section*{Introduction}
	
	A deeply rooted intuition about smooth manifolds is that of spaces that become linear spaces in the infinitesimal neighbourhood of each point. On the infinitesimal scale the geometry underlying a manifold is thus affine geometry.
	To make this intuition precise requires a good theory of infinitesimals as well as defining precisely what it means for two points on a manifold to be infinitesimally close. As regards infinitesimals we make use of \emph{Synthetic Differential Geometry} (SDG) and adopt the neighbourhoods of the diagonal from Algebraic Geometry to define when two points are infinitesimally close. The key observations on how to proceed have been made by Kock in \cite{SGM}: 1)~The first neighbourhood of the diagonal exists on formal manifolds and can be understood as a symmetric, reflexive relation on points, saying when two points are infinitesimal neighbours, and 2)~we can form affine combinations of points that are mutual neighbours.
	
	It remains to make precise in which sense a manifold becomes a model of the theory of affine spaces. This has been done in \cite{IMAT}. Firstly, one abstracts from Kock's infinitesimal simplices of mutual infinitesimally neighbouring points to what is called an \emph{infinitesimal structure}. (See also section~\ref{sec:i-aff-space} for a definition.) An infinitesimal structure serves then as the domain of definition for the operations of affine combinations. A space together with an infinitesimal structure (i-structure) and an action of the clone of affine operations on that infinitesimal structure is called an \emph{infinitesimally affine space} (i-affine space).  
	
	Formal manifolds and affine schemes (considered as either duals of commutative rings, or $\Coo$-rings) are examples of i-affine spaces. The i-structures are generated by the first neighbourhood of the diagonal. In this paper we shall construct i-structures from the $k$th-order neighbourhoods of the diagonal on $R^n$ for a ring $R$ satisfying the Kock-Lawvere axioms for higher-order infinitesimals. The definition of these i-structures are guided by the requirements that these i-structures are preserved by all maps $f:R^n\to R^m$ (hence can be defined on formal manifolds as well) and that the affine structure of $R^n$ restricts to an i-affine space on each higher-order i-structure. Both of these hold true for the i-structure generated by the first neighbourhood of the diagonal. In contrast to the first neighbourhood of the diagonal the i-affine structures on the higher-order neighbourhoods are not preserved by all maps anymore. Therefore, whereas a manifold carries all the higher-order i-structures, an i-affine structure has to be imposed as an additional piece of data.
	
	We show that any second-order i-affine structure on a manifold induces a symmetric affine connection, and, conversely, any symmetric affine connection extends to a second-order i-affine structure in such a way that the latter is of the same affine-algebraic form as the canonical connection on an affine space. Furthermore, as we are dealing with affine connections on points, we shall also discuss existence results of such connections, and hence the existence of second-order i-affine structures on (smooth) manifolds.
	
	
	\section{Infinitesimally affine spaces}\label{sec:i-aff-space}
	
	We shall work mostly within naive axiomatic SDG, as it is done in \cite{SGM}, for example. Let $A$ be a space. An \emph{i-structure} on $A$ amounts to give an $n$-ary relation $A\la n\ra$ for each $n\in\IN$ that defines which $n$ points in $A$ are considered as being `infinitesimally close' to each other.
	\begin{definition}[i-structure]
		Let $A$ be a space. An \emph{\textbf{i-structure}} on $A$ is an $\IN$-indexed family $n\mapsto A\la n\ra\subseteq A^n$ such that
		\begin{enumerate}[(1)]
			\item $A\la 1\ra = A$, $A\la 0\ra=A^0=1$ (the `one point' space, or terminal object)
			\item For every map $h:m\to n$ of finite sets 
			and every $(P_1,\ldots,P_n)\in A\la n\ra$ we have $ (P_{h(1)},\ldots,P_{h(m)})\in A\la m\ra$ 
		\end{enumerate}
	\end{definition}
	The first condition is a normalisation condition. The second condition  makes sure that the relations are compatible: if we have a family of points that are infinitesimally close to each other, then so is any subfamily of these points, or any family created from repetitions. In particular, we obtain that the $A\la n\ra$ are symmetric and reflexive relations. An $n$-tuple $(P_1,\ldots,P_n)\in A^n$ that lies in $A\la n\ra$ will be denoted by $\la P_1,\ldots,P_n\ra$ and we shall refer to these as \emph{i-$n$-tuples}. A map $f:A\to X$ that maps i-$n$-tuples to i-$n$-tuples for each $n\in\IN$, i.e. $f^n(A\la n\ra)\subseteq X\la n\ra$, is called an \emph{i-morphism}.
	
	Two trivial examples of i-structures on $A$ are the discrete and indiscrete i-structure obtained by taking $A\la n\ra$ to be the diagonal $\Delta_n$, respectively the whole $A^n$. The i-structures that are of main interest in SDG are the i-structures generated by the first neighbourhood of the diagonal (as relations). We call them \emph{nil-square i-structures}. For example, let $R$ be a ring\footnote{All rings are assumed to be commutative.}. Recall that
	$$D(n)=\set{(d_1,\ldots,d_n)\in R^n}{d_id_j=0,\ 1\leq i,j\leq n}$$
	On $R^n$ the first neighbourhood of the diagonal is given by 
	$$
	\set{(P_1,P_2)}{P_2-P_1\in D(n)}
	$$ 
	This is a symmetric and reflexive relation and we can construct an i-structure from it: take the first neighbourhood of the diagonal as $R^n\la 2\ra$ and define the \emph{nil-square i-structure} on $R^n$ by
	$$R^n\la m\ra=\set{(P_1,\ldots,P_m)}{(P_i,P_j)\in R^n\la 2\ra,\ 1\leq i,j \leq m}$$    
	This i-structure is thus \emph{generated} by $R^n\la 2\ra$. Not all i-structures $A\la-\ra$ of interest need to be generated by $A\la 2\ra$. We will see such examples in section~\ref{sec:higher-i-structure}.
	
	If the ring $R$ satisfies the \emph{Kock-Lawvere axiom}, that is for every map $t:D(n)\to R$ there are unique $a_0,\ldots,a_n\in R$ such that $$t(d_1,\ldots,d_n)=a_0+\sum_{j=1}^n a_j d_j,\qquad (d_1,\ldots,d_n)\in D(n),$$ 
	then every map $f:R^n\to R^m$ is an i-morphism of the nil-square i-structures. This is due to the following two facts: linear maps $R^n\to R^m$ map $D(n)$ to $D(m)$, and for $P_2-P_1\in D(n)$ 
	\begin{equation}\label{eq:differential}
	    f(P_2)-f(P_1)=\partial f(P_1)[P_2-P_1] 
	\end{equation}
	where $\partial f(P_1)$ denotes the \emph{derivative of $f$} at $P_1$. The stated property of linear maps can be checked by direct computation; the existence and uniqueness of the linear map $\partial f(P_1)$ are both a consequence of the Kock-Lawvere axiom.
	%
	
	The nil-square i-structure induces i-structures on subspaces $U\into R^n$ by restriction. For \emph{formally open subspaces} $U\into R^n$, which are stable under infinitesimal perturbations at each point (see \cite[I.17]{SDG} or \cite[def.~3.2.5]{IMAT} for a definition), each map $f:U\to R^m$ has a derivative; hence every map $f:U\to V$ between formally open subspaces is an i-morphism. Furthermore, it is possible to glue the i-structures on formally open subspaces together to get an i-structure on a formal manifold and show that every map between formal manifolds is an i-morphism. (See \cite[prop.~I.17.5]{SDG} and \cite[thm.~3.2.8]{IMAT} for proofs.)  
	\begin{definition}[i-affine space]
		Let $A\la-\ra$ be an \emph{i-structure} on $A$. Set $\calA(n) = \set{(\lambda_1, . . . , \lambda_n)\in R^n}{\sum_{j=1}^{n}\lambda_j=1}$. The space $A$ is said to be an \emph{\textbf{i-affine space}} (over $R$), if for every $n\in\IN$ there are operations
		$$\calA(n)\times A\la n\ra\to A,\qquad ((\lambda_1,\ldots,\lambda_n),\la P_1,\ldots,P_n\ra)\mapsto \sum_{j=1}^{n} \lambda_j P_j$$
		satisfying the axioms
		\begin{itemize}
			\item (\emph{Neighbourhood}) Let $\lambda^k\in \calA(n)$, $1\leq k\leq m$. If $\la P_1,\ldots,P_n \ra\in A \la n\ra$  then
			$$ \big(\sum_{j=1}^n \lambda^1_j P_j,\ldots, \sum_{j=1}^n \lambda^m_j P_j\bigr)\in A\la m\ra$$
			\item (\emph{Associativity}) Let $\lambda^k\in \calA(n)$, $1\leq k\leq m$, $\mu \in\calA(m)$ and $\la P_1,\ldots,P_n \ra\in A \la n\ra$. 
			$$ \sum_{k=1}^m\mu_k\big(\sum_{j=1}^n \lambda^k_j P_j\bigr)=\sum_{j=1}^n \bigl(\sum_{k=1}^m\mu_k\lambda^k_j\bigr)P_j$$
			(Note that the left-hand side is well-defined due to the neighbourhood axiom.)
			\item (\emph{Projection}) Let $n\geq 1$ and let $e^n_k\in R^n$ denote the $k$th standard basis vector for $1\leq k\leq n$. For every  $\la P_1,\ldots,P_n \ra\in A \la n\ra$ it holds
			$$\sum_{j=1}^n (e^n_k)_j P_j = P_k$$
			In particular, we have for $n=1$ that $1P=P$, $P\in A$.   
		\end{itemize}
	\end{definition}
	The neighbourhood axiom makes sure that we can compose affine combinations as we are used to, provided we are working over a fixed i-tuple. The associativity and projection axioms make sure the algebra of affine combinations follows the same rules as in all the $R^n$. A consequence of the neighbourhood axiom is that every i-tuple generates an affine space over $R$. This makes precise the statement that the geometry of the space $A$ is affine on the infinitesimal scale.
	
	It is not difficult to show by direct calculation that the affine space $R^n$ satisfies the neighbourhood axiom for the nil-square i-structure making it an i-affine space\footnote{This is also a consequence of the more general \cite[cor.~3.1.6 and 2.3.3]{IMAT}.}. Moreover, due to (\ref{eq:differential}) it follows that every map $f:R^n\to R^m$ preserves not only the nil-square i-structure but the i-affine combinations as well. Each map $f$ is an \emph{i-affine map}.
	
	The i-affine structure of $R^n$ restricts to its formally open subspaces. Due to (\ref{eq:differential}) all maps between formally open subspaces become i-affine maps for these i-structures. Like with the i-structures also the i-affine structures on formally open subspaces can be glued together to an i-affine structure on a formal manifold. All maps between formal manifolds become i-affine maps for these i-affine structures \cite[thm.~3.2.8]{IMAT}. Any manifold in the sense of classical differential geometry is a formal manifold\footnote{This is to be understood in the context of well-adapted models of SDG \cite{WAM}, where we have a fully faithful embedding of the category of smooth manifolds into a Grothendieck topos that admits a model of the Kock-Lawvere axioms. This embedding maps the real line $\IR$ to $R$, analytical derivatives to derivatives in SDG and it maps open covers to covers by formally open spaces \cite[III.3]{SDG}.}, so any manifold is an i-affine space and any smooth map between manifolds is i-affine.
	
	Affine schemes (considered as either duals of commutative rings, or $\Coo$-rings) become examples of i-affine spaces over their respective nil-square i-structure \cite[cor.~2.3.3 and 3.1.6]{IMAT}. Every morphism of affine schemes becomes an i-morphism. Affine $\Coo$-schemes, for example, form a category of spaces generalising smooth manifolds. Besides manifolds the category fully faithfully embeds locally closed subsets of Euclidean space with smooth maps between them \cite[prop.~1.5]{SIA}. This provides us with a wealth of examples of i-affine spaces. Furthermore, i-affine spaces are surprisingly well-behaved under taking colimits of the underlying spaces \cite[chap.~2.6]{IMAT}, \cite{IMAT_paper}. This and their algebraic nature makes them a suitable type of space to study geometric notions based on infinitesimals.   
	
	\section{Higher-order infinitesimal structures}\label{sec:higher-i-structure}
	
	The important examples of i-structures so far have all been the nil-square i-structures, which are constructed from the first neighbourhood of the diagonal. In this section we wish to define i-structures $A_k=A_k\la-\ra$ on $A=R^n$ such that $A_k\la 2\ra$ is the $k$th neighbourhood of the diagonal
	$$
	\set{(P_1,P_2)}{P_2-P_1\in D_k(n)}
	$$ 
	where $D_k(n)$ is the space of $k$th-order infinitesimals     
	$$D_k(n)=\set{(d_1,\ldots,d_n)\in R^n}{\text{any product of $(k+1)$ $d_j$ vanishes}}$$
	The i-structures $A_k\la-\ra$ shall satisfy
	\begin{enumerate}[1)]
		\item All maps $f:R^n\to R^m$ become i-morphisms for the respective $k$th-order i-structures on $R^n$ and $R^m$
		\item The affine space $A=R^n$ becomes an i-affine space over $A_k\la-\ra$.
	\end{enumerate}
	To be able to study 1) we assume henceforth that $R$ is a $\IQ$-algebra that satisfies the Kock-Lawvere axiom for all the $D_k(n)$, $k,n\geq 1$. This amounts to say that each map $t:D_k(n)\to R$ is a polynomial function for a uniquely determined polynomial in $R[X_1,\ldots,X_n]$ of total degree $\leq k$. An important consequence is that every map $f:A\to R^m$ has a \emph{Taylor representation} 
	$$ f(P)-f(Q)= \sum_{\ell=1}^{k}\frac{1}{\ell!}\partial^\ell f(Q)[P-Q]^\ell
	$$
	for $P-Q\in D_k(n)$. Here $\partial^\ell f(Q)$ stands for the $\ell$th derivative of $f$ at $Q$, which is a symmetric $\ell$-linear map $(R^n)^\ell\to R^m$. Writing $\phi[v]^\ell$ for an $\ell$-linear map $\phi$ means that we evaluate it on the $\ell$-tuple $(v,\ldots,v)$. The following characterisation of $D_k(n)$ in \cite[prop.~1.2.2]{SGM} will be useful
	$$
		D_k(n)=\set{d\in R^n}{\phi[d]^{(k+1)}=0\text{ for all $(k+1)$-linear $\phi:(R^n)^{(k+1)}\to R$}}
	$$
	Let $V\cong R^n$ and $k\geq 1$. We define $\DN_k(V)$ to be the space
	\begin{multline*}
		\DN_k(V)=\{(v_1,\ldots, v_{k+1}) \in \D_k(V)^{(k+1)}\,\mid\\ \text{For any $(k+1)$-linear map }\phi: V^{(k+1)} \to R,\ \phi(v_1,\ldots,v_{k+1})=0\}
	\end{multline*}
	In the subsequent definition we will use $A=R^n$ to mean the (affine) space $R^n$ and $V=R^n$ to mean the $R$-vector space $R^n$.
	\begin{definition}[$k$th-order i-structure on $R^n$]\label{def:kth-i-str} Let $A=V=R^n$ and $k\geq 1$. We define the \emph{\textbf{$k$th-order i-structure}} $A_k$ on $A$ by
		\begin{enumerate}[(1)]
			\item $A_k\la 1\ra = A$, $A_k\la 0\ra = A^0=1$
			\item For $m\geq 2$
			\begin{multline*}
			A_k\la m \ra = \{(P_1,\ldots, P_m) \in A^m \mid (P_{i_1} - P_{j_1},\ldots, P_{i_{k+1}}-P_{j_{k+1}})\in \DN_k(V),\\ \text{ for all } i_\ell, j_\ell \in \{1,\ldots, m\},\ i_\ell\neq j_\ell, 1\leq \ell\leq k+1 \}
			\end{multline*}
		\end{enumerate}
	\end{definition}
	From the definition it follows readily that each $A_k$ is indeed an i-structure and that  
	$$A_k\la 2\ra = \set{(P_1,P_2)\in A^2}{P_2-P_1\in D_k(n)}$$
	is the $k$th neighbourhood of the diagonal as desired. 
	%
	%
	
	Note that the first-order i-structure $A_1$ is smaller than the nil-square i-structure on $A=R^n$ for $n>1$, i.e. $A_1\la m\ra\subseteq A\la m\ra$ for all $m\in\IN$. Indeed, both i-structures agree up to $m=2$, but $\la P_1,P_2,P_3\ra \in A\la 3\ra$ if and only if $\phi(P_i-P_k,P_j-P_k)=0$ for every \emph{symmetric} bilinear form $\phi$ and every $1\leq i,j,k\leq 3$ \cite[prop.~1.2.12]{SGM}. We would have had both i-structures agree, if we had restricted to symmetric $(k+1)$-linear forms in the definition of $\DN_k(V)$. The reason for not doing so is that this i-structure is not provably preserved by all maps $f:R^n\to R^m$ for $k\geq 2$.
	
	\begin{theorem}\label{thm:kth-i-morph} Every map $f:R^n\to R^m$ is an i-morphism for the respective $k$th-order i-structures.  
	\end{theorem}
	\begin{proof} 
		To avoid any more overload of notation with indices we will denote the $k$th-order i-structure on $R^n$ with $A_k$ and the one on $R^m$ with $B_k$. Moreover, we set $V_A=R^n$ and $V_B=R^m$. Let $\la P_1,\ldots,P_m\ra \in A_k\la m\ra$ for an index $m\geq 2$. We have to show 
		$$
			\la f(P_1),\ldots,f(P_m)\ra \in B_k\la m\ra
		$$
		By definition this amounts to show
		$$
			\phi(f(P_{i_1}) - f(P_{j_1}),\ldots, f(P_{i_{k+1}})-f(P_{j_{k+1}}))=0
		$$
		for all $i_\ell, j_\ell \in \{1,\ldots, m\}$,  $1\leq \ell\leq k+1$ and any $(k+1)$-linear form $\phi$ on $V_B$. Since each $P_{i_\ell} - P_{j_\ell}\in D_k(n)$ we can apply Taylor expansion
		$$
			f(P_{i_\ell})-f(P_{j_\ell})= \sum_{j=1}^{k}\frac{1}{j!}\partial^j f(Q)[P_{i_\ell}-P_{j_\ell}]^j
		$$
		Substituting each $f(P_{i_\ell})-f(P_{j_\ell})$ with its respective Taylor expansion in $\phi$ and applying multilinearity to expand the $k+1$ sums yields a sum of multilinear forms on $V_A$ of the order $(k+1)$ or higher with arguments being combinations of $P_{i_\ell} - P_{j_\ell}$ for $i_\ell, j_\ell \in \{1,\ldots, m\}$,  $1\leq \ell\leq k+1$. Because of $\la P_1,\ldots,P_m\ra \in A_k\la m\ra$ each such multilinear form evaluates to $0$, hence does the sum. This shows that 
		$$
			\phi(f(P_{i_1}) - f(P_{j_1}),\ldots, f(P_{i_{k+1}})-f(P_{j_{k+1}}))=0
		$$
		as required. We conclude that $f$ is an i-morphism as claimed.
	\end{proof}
	The proof of the preceding proposition clarifies why we need to define $\DN_k(V)$ using $(k+1)$-multilinear maps and not just the symmetric ones: even though each multilinear map in the Taylor expansion of $f(P_{i_\ell})-f(P_{j_\ell})$ is symmetric, the expansion is a sum over multilinear maps of different degrees. Once we expand
	$$\phi(f(P_{i_1}) - f(P_{j_1}),\ldots, f(P_{i_{k+1}})-f(P_{j_{k+1}}))=0$$
	into a sum of multilinear maps, those multilinear maps will be compositions of $\phi$ with multilinear maps of different degrees and hence \emph{not symmetric anymore}, in general. For example, consider $k=2$, $\la P_1,P_2,P_3 \ra \in A_2\la 3\ra$ and a symmetric trilinear form $\phi$ on $V_B$. After a tedious but straight-forward calculation one obtains
	\begin{multline*}
	\phi(f(P_2) - f(P_1),f(P_3) - f(P_1), f(P_3) - f(P_2))= \\ \frac{1}{2}\bigl(\phi(\partial f(P_1)[P_2-P_1],\partial f(P_1)[P_2-P_1],\partial^2f(P_1)[P_3-P_1]^2)\\
	- \phi(\partial f(P_1)[P_3-P_1],\partial f(P_1)[P_3-P_1],\partial^2f(P_1)[P_2-P_1]^2)\bigr)
	\end{multline*}
	The right hand side is not provably equal to $0$ for all symmetric trilinear forms $\phi$, in general. Therefore, defining $\DN_2(V_A)$ using symmetric trilinear forms only instead of all trilinear forms would make us unable to prove that all $f$ preserve the 2nd-order i-structure, for example.
	%
	%
	\begin{theorem}\label{thm:kth-i-aff} The affine structure on $A=R^n$ restricts to the $k$th-order i-structure $A_k$, making $A_k$ an i-affine subspace of the affine space $A$ (equipped with the indiscrete i-structure).    
	\end{theorem}
	\begin{proof}
		We shall make use of the notation from the proof of the preceding proposition. To show $A_k$ an i-affine subspace of $A$ it suffices to show that the affine operations on $A$ satisfy the neighbourhood axiom for $A_k$.
		
		Let $\lambda^i\in \calA(n)$ for $1\leq i\leq m$ and $\la P_1,\ldots,P_n \ra\in A_k \la n\ra$. We have to show
		$$ \big\la\sum_{j=1}^n \lambda^1_j P_j,\ldots, \sum_{j=1}^n \lambda^m_j P_j\bigr\ra\in A_k\la m\ra$$
		Let $\phi$ be a $(k+1)$-linear form on $V_A$ and $i_\ell, j_\ell \in \{1,\ldots, m\}$ for all $1\leq \ell\leq k+1$. Using $\sum_{j=1}^n \lambda^{i}_j =1$, $1\leq i\leq m$  yields
		\begin{align*}
		&\phi\bigl(\sum_{i=1}^n \lambda^{i_1}_i P_i-\sum_{j=1}^n \lambda^{j_1}_j P_j,\ldots, \sum_{i=1}^n \lambda^{i_{k+1}}_i P_i-\sum_{j=1}^n \lambda^{j_{k+1}}_j P_j\bigr)\\
		&=\phi\bigl(\sum_{i,j=1}^n \lambda^{i_1}_i\lambda^{j_1}_j (P_i- P_j),\ldots, \sum_{i,j=1}^n \lambda^{i_{k+1}}_i \lambda^{j_{k+1}}_j (P_i-P_j)\bigr)
		\end{align*}
		Applying the multilinearity of $\phi$ yields a sum of $(k+1)$-linear forms with arguments being combinations of $P_{i_\ell} - P_{j_\ell}$ for $i_\ell, j_\ell \in \{1,\ldots, n\}$,  $1\leq \ell\leq k+1$, which all evaluate to zero by assumption. We conclude
		$$ \big\la\sum_{j=1}^n \lambda^1_j P_j,\ldots, \sum_{j=1}^n \lambda^m_j P_j\bigr\ra\in A_k\la m\ra$$
		as required.
	\end{proof}
	%
	%
	The definitions of the $k$th-order i-structure $A_k$ together with theorems~\ref{thm:kth-i-morph} and \ref{thm:kth-i-aff} can be generalised to a formally open subspace $A$ of $R^n$, directly. This allows us to glue together the $k$th-order i-structures to a $k$th-order i-structure on a formal manifold and all maps between formal manifolds will preserve that structure. 
	
	\begin{theorem} \label{thm:kth-i-struc} Let $A$ be a formal manifold and $k\geq 1$.
		\begin{enumerate}[(i)]
			\item $A$ carries a unique i-structure $A_k$ with the universal property that any map $f:A\to M$ to a space $M$ equipped with an i-structure is an i-morphism $f:A_k\to M$ if and only if for every formally open subspace $\iota: U\into A$ that is also a formally open subspace of $R^n$ (i.e. a \emph{chart} of $A$) the restriction of $f$ along $\iota$ is an i-morphism $U_k\to M$.
			
			Here $U_k$ denotes the $k$th-order i-structure on $U$ as a formally open subspace of $R^n$; i.e. the pullback of the $k$th-order i-structure of $R^n$ to $U$.  
			%
			\item All maps between formal manifolds become i-morphisms for the respective $k$th-order i-structures. 
		\end{enumerate}
	\end{theorem}
	\begin{proof}
		\begin{enumerate}[(i)]
			\item (Essentially, this part is theorem~2.6.19 in \cite{IMAT} applied to the i-structure only.) For each $n\geq 1$ we define $A_k\la n \ra$ as the join of the images of $U_k\la n\ra$ for each chart $\iota: U\into A$ of $A$. It is easy to see that this yields an i-structure on $A$ with the desired universal property.    
			\item Let $f: A\to M$ be a map between two formal manifolds equipped with the $k$th-order i-structure as defined in (i) and $\la P_1,\ldots,P_n\ra\in A_k\la n\ra$. By construction there is an $A$-chart $\iota: U\into A$, $\phi: U\into R^n$, and $\la x_1,\ldots,x_n\ra\in U_k\la n\ra$ such that $\iota(x_\ell)=P_\ell$, $1\leq \ell\leq n$.
			
			We also find an $M$-chart $j:V\into M$ containing $f(P_1)$. Pulling back $j$ along $f$ yields a formally open subspace $f^*j:f^{-1}(V)\into M$, which becomes a chart after taking the intersection with $\iota$
			$$
				\iota^*f^*j: U\cap f^{-1}(V)\into A, \quad (\iota^*f^*j)^*\phi: U\cap f^{-1}(V)\into R^n
			$$
			(Recall that formally open subspaces are stable under pullback.) Let $W = U\cap f^{-1}(V)$. The restriction of $f: W\to V$ is a map between formally open subspaces of $R^n$ and $R^m$, respectively, and thus an i-morphism by theorem~\ref{thm:kth-i-morph} and the constructions of $W_k$ and $V_k$. Since $x_1\in W\subset U$ and $W$ is a formally open subspace of $U$, we find $\la x_1,\ldots,x_n\ra\in W_k\la n\ra$ and hence $\la f(x_1),\ldots,f(x_n)\ra \in V_k\la n\ra$; but this implies that 
			$$
				\la f(P_1),\ldots,f(P_n)\ra = \la j(f(x_1)),\ldots,j(f(x_n))\ra \in M_k\la n\ra
			$$    
			and that $f$ is an i-morphism as claimed.  
		\end{enumerate} 
	\end{proof}

	\begin{remark}\label{rem:kth-i-struc} Part~(i) of the preceding theorem states in simpler terms that $f:A\to M$ is an i-morphism, if and only if it is an i-morphism on the charts. Instead of forming the union over all charts, $A_k$ can be also defined as the union over a covering family, i.e. an \emph{atlas}. Moreover, $f$ is an i-morphism if and only if all its restriction to the charts of the atlas are i-morphisms. 
	
	Indeed, any chart of $\iota: U\into A$ can be covered by restrictions of charts of the chosen atlas, which are formally open subspaces of both $A$ and some $R^n$. The same argument as presented in the proof of (ii) above shows that $\iota$ is an i-morphism when applied to $U$ and charts of the atlas.  
	\end{remark}

	However, note that theorem~\ref{thm:kth-i-struc} does not extend to the i-affine structures, i.e. maps are not going to preserve the i-affine structure on $U_k$ for a formally open subspace $U\into R^n$, in general. Only special classes of maps will have that property and these classes will depend on $k$. Indeed, for $k\geq 2$ the Taylor expansion of $f$ contains quadratic terms and higher, hence can only preserve affine combinations up to quadratic and higher-order terms. Therefore, unlike $R^n$ a formal manifold does \emph{not} carry a canonical i-affine structure on its canonical $k$th-order i-structure.    
	
	Let $A=R^n$ or, more generally, a formally open subspace of $R^n$. Besides the i-affine structure over the nil-square i-structure we have now i-affine structures over each $k$th-order i-structure. It is readily seen from the definitions that $D_k(n)\subseteq D_{k+1}(n)$ and $A_k\la m\ra\subseteq A_{k+1}\la m\ra$. The identity map $1_A: A\to A$ thus induces an i-affine embedding $A_k\into A_{k+1}$. If $A$ is a formal manifold, then this embedding remains an i-morphism.
	
	\begin{corollary} Let $A_k$ denote the $k$th-order i-structure on a formal manifold $A$, $k\geq 1$.
		\begin{enumerate}[(i)]
			\item The identity map $1_A:A\to A$ induces i-embeddings 
			$A_k\into A_{k+1}$
			\item In the case of $A$ being a formally open subspace of $R^n$ the inclusions $A_k\into A_{k+1}$ become i-affine maps for the i-affine structures on each $A_k$. 
		\end{enumerate}
	\end{corollary}
	
	\begin{remark}\label{rem:nil-square}
		As regards the nil-square structure on $A=R^n$, depending on the dimension of $A$ it might not be provably contained in any of the $A_k$. Indeed, we find that 
		$$A\la m\ra \subseteq A_{m-1}\la m\ra$$
		This follows from the fact that for any $m$-tuple of points you can only form $m-1$ different difference vectors. Hence any argument for an $m$-linear form will contain at least one repetition of a difference vector, and thus has to vanish.
		
		Is $m-1$ a strict bound for the inclusion of $A\la m\ra$? We analyse the behaviour of multilinear forms on the nil-square i-structure more carefully. Recall that $\la P_1,P_2,P_3\ra\in A\la 3\ra$ if and only if $\phi[u,v]=0$ for any symmetric bilinear form $\phi$, where $u=P_{j_1}-P_{i_1}$ and $v=P_{j_2}-P_{i_2}$. For a general bilinear form this implies that $\phi[u,v]=-\phi[v,u]$. Therefore, if $\la P_1,\ldots,P_{m+1}\ra\in A\la m+1\ra$ and $(v_1,\ldots,v_{m})$ is an $m$-tuple of vectors with $v_\ell=P_{i_\ell}-P_{j_\ell}$ for some $1\leq i_\ell,j_\ell\leq m+1$, then any $m$-form $\phi$ is alternating on $(v_1,\ldots,v_m)$. This means that as long as we can find $m+1$ points $\la P_1,\ldots,P_{m+1}\ra$, which difference vectors have a determinant that is not provably equal to $0$, we can find an $m$-linear form that does not provably evaluate to zero on the difference vectors showing that $A\la m+1\ra$ is not contained in $A_{m}\la m+1\ra$ provided $m\leq n$.
		
		For each $d_1\in D=D(1)$ and any $m\geq 1$ the Kock-Lawvere axiom guarantees the existence of $d_2,\ldots,d_m\in D$ such that their product $d_1\cdots d_m$ is not provably equal to zero. We construct $m$ vectors $v_j\in R^m$
		$$
		v_1= d_1\, e_1,\ 
		v_2 = d_2\, e_2,\ \ldots,\ 
		v_{m}=d_{m}\,e_m,
		$$
		where the $e_j$ form the \emph{standard basis} of $R^m$. The determinant $\det[v_1,\ldots,v_m]$ evaluates to $d_1\cdots d_m$. 
		
		Let $n$ be the dimension of $A$. Suppose $m\leq n$, then by extending the components of each $v_j$ with $n-m$ zeros we obtain the desired $m+1$ points $\la 0,v_1,\ldots,v_{m}\ra\in A\la m+1\ra$. Pulling back the determinant along the projection $R^n\to R^m$ onto the first $m$ components yields an $m$-linear form that does not provably evaluate to zero on the difference vectors $(v_1,\ldots,v_m)$ as claimed.
		%
		%
	\end{remark}
	
	\section{Affine connections and 2nd-order i-affine structures}\label{sec:connections}
	
	In differential geometry affine connections on a manifold come in three equivalent notions: a geometric notion of parallel transport of tangent vectors along paths, and two algebraic notions; that of a covariant derivative on vector fields and the horizontal subbundle of the iterated tangent bundle. In SDG we can study these notions from the infinitesimal viewpoint. A tangent vector at a point $P$ is an `infinitesimal piece' of a path:
	$t:D\to A$
	with $t(0)=P$. Geometrically, a parallel transport of a tangent vector $t_1$ along a path $\gamma: [0,1]\to A$ amounts to an `infinitesimal thickening' of $\gamma$ in the direction of $t_1$, that is a map
	$$P_\gamma(t_1): D\times [0,1]\to A $$
	If we replace $\gamma$ with a tangent vector $t_2$ over the same base point as $t_1$ the situation becomes symmetric
	$$P_{t_1}(t_2): D\times D\to A$$
	From the infinitesimal viewpoint an affine connection is thus essentially a mapping $\nabla$ that takes a pair of tangent vectors $(t_1,t_2)$ over the same base point and assigns them a tangent square $\nabla(t_1,t_2)=P_{t_1}(t_2)$ over that base point such that the principal axes of this tangent square are $t_1$ and $t_2$. By noting that the iterated tangent bundle $TTA\to A$ is the bundle of tangent squares $A^{D\times D}\to A$ one can readily relate the affine connection with a covariant derivative and the horizontal subbundle \cite[chap.~5]{Lav}, \cite[chap.~4.6]{SGM}.
	
	For a formal manifold $A$ the points are geometrically more fundamental than tangent vectors. Indeed, one can show that the vector space structure on each tangent space $T_PA$ is a pointwise linear structure on the maps $D\to A$ derived from $A$ being infinitesimally linear at $P$ \cite[chap.~4.2]{SGM}, \cite[chap.~3.3.2]{IMAT}. Like an affine connection completes two tangent vectors to a tangent square, an affine connection for points takes three points $P,Q,S$ and completes them to a parallelogram $PQRS$ \cite[chap.~2.3]{SGM}. Here $\la P,Q\ra$ and $\la P,S\ra$ are first-order neighbours, but $Q$ and $S$ don't need to be. The resulting point $R$ is a first-order neighbour of $P$ and of $Q$, hence it is a second-order neighbour of $P$. If we follow \cite{SGM} and denote the point $R$ by $\lambda(P,Q,S)$ then an \emph{affine connection} (on points) $\lambda$ is a map mapping triples $(P,Q,S)$ with $\la P,Q\ra, \la P,S\ra\in A\la 2\ra$ to a point $\lambda(P,Q,S)$ such that
	\begin{align*}
		\lambda(P,Q,P) &= Q\\
		\lambda(P,P,S) &= S\\
	\end{align*}
	These properties are sufficient to derive the other nil-square neighbourhood relationships \cite[chap.~2.3]{SGM}. An affine connection is called \emph{symmetric}, if 
	$$\lambda(P,Q,S)=\lambda(P,S,Q) $$
	For $A=R^n$ a symmetric affine connection is induced by its affine structure 
	$$\lambda(P,Q,S)=Q+S-P$$
	Geometrically, this is the addition of vectors using parallel transport to construct a vector parallelogram at $P$. In fact, any i-affine structure on $A_2$ induces a symmetric affine connection in this way.
	
	\begin{proposition}\label{prop:2nd-iaff-connection} Let $A$ be a formal manifold that admits an i-affine structure on $A_2$, then $A$ admits a symmetric affine connection on points.
	\end{proposition}
	\begin{proof}
		We wish to define the symmetric affine connection $\lambda$ by
		$$\lambda(P,Q,S):=Q+S-P$$
		where the right hand side denotes the i-affine combination in $A_2$. For this to be well-defined we need to show $\la P,Q,S\ra\in A_2\la 3\ra$. We work in a chart. First note that $Q-P,S-P,Q-S\in D_2(n)$. Let $\phi$ be a trilinear map. We find
		$$\phi[Q-P,S-P,Q-S] = \phi[Q-P,S-P,Q-P]-\phi[Q-P,S-P,S-P]=0$$
		as the two trilinear maps on the right hand side are quadratic in $Q-P\in D(n)$, respectively in $S-P\in D(n)$. This is sufficient to show $\la P,Q,S\ra\in A_2\la 3\ra$.  The defining properties showing $\lambda$ an affine connection are immediate consequence of the algebra of affine combinations. 
	\end{proof}
	
	We wish to show the converse, i.e. that any symmetric affine connection $\lambda$ on a formal manifold $A$ extends to a 2nd-order i-affine structure. To show this we shall proceed in two steps. First we show that this holds on a formally open subspace of $U\subseteq R^n$. Then we show that for any formally open subspace $V\subseteq R^n$ with an embedding $\iota: V\into U$ the 2nd-order i-affine structure defined on $V$ by $\lambda$ is preserved by $\iota$. This allows us to glue the 2nd-order i-affine structures together to a 2nd-order i-affine structure on the formal manifold $A$ (theorem~3.2.8\footnote{Although theorem~3.2.8 refers to the nil-square i-structure only, due to being formally open and theorem~\ref{thm:kth-i-struc} all the required properties of charts also hold for the 2nd-order i-affine structure. The assertion of theorem~3.2.8 can thus be extended to the 2nd-order i-affine structure when combining the original proof with the subsequent lemmas.} or theorem~2.6.19 in \cite{IMAT}).
	
	Let $\lambda$ be a connection on $U$. It is not difficult to show that
	$$\lambda(P,Q,S)=Q+S-P + \Gamma_P[Q-P,S-P]$$
	for a symmetric bilinear map $\Gamma_P$ \cite[chapter~2.3]{SGM}, which we will refer to as \emph{Christoffel symbols} of the connection as it is done in \cite{SGM}. For each $n\geq 1$ we define an action of $\calA(n)$ on $U_2\la n\ra$ by
	$$\lambda\cdot\la P_1,\ldots,P_n\ra = \sum_{j=1}^n\lambda_jP_j + \frac{1}{2}\bigl(\Gamma_{P_1}\bigl[\sum_{j=1}^n\lambda_jP_j-P_1\bigr]^2-\sum_{j=1}^n\lambda_j
	\Gamma_{P_1}[P_j-P_1]^2\bigr)$$
	Firstly, note that due to $\la P_1,\ldots,P_n\ra\in U_2\la n\ra$ and
	$$\sum_{j=1}^n \lambda_j P_j=\bigl(1-\sum_{j=2}^{n} \lambda_j\bigr)P_1+\sum_{j=2}^{n} \lambda_j P_j=P_1+\sum_{j=2}^{n} \lambda_j (P_j-P_1)$$
	the vector
	$$\lambda\cdot\la P_1,\ldots,P_n\ra - P_1= \sum_{j=2}^n\lambda_j(P_j-P_1) + \frac{1}{2}\bigl(\Gamma_{P_1}\bigl[\sum_{j=2}^n\lambda_j(P_j-P_1)\bigr]^2-\sum_{j=1}^n\lambda_j
	\Gamma_{P_1}[P_j-P_1]^2\bigr)$$
	lies in $D_2(n)$. Furthermore, for any $\lambda^1,\ldots,\lambda^m\in \calA(n)$ we have
	$$\la \lambda^1\cdot\la P_1,\ldots,P_n\ra,\ldots,\lambda^m\cdot\la P_1,\ldots,P_n\ra\ra \in U_2\la m\ra$$
	which shows the neighbourhood axiom. For all the standard  
	basis vectors $e^n_k\in \calA(n)$ we find
	$$\Gamma_{P_1}\bigl[\sum_{j=1}^n(e^n_k)_jP_j-P_1\bigr]^2-\sum_j(e^n_k)
	\Gamma_{P_1}[P_j-P_1]^2=0$$
	so the projection axiom holds true as well. The proof of the associativity axiom involves a longer calculation, and we will give only the most important steps. The main techniques used in this calculation are Taylor-expansion and multilinear algebra  of nil-potents we have been using a lot already. We exhibit these types of arguments in more detail while showing that these actions by affine combinations on formally open subsets are compatible first, as the calculations are simpler than in the proof of associativity.
	
	\begin{lemma} Let $U$, $V$ be formally open subsets of $R^n$, $\iota: V\into U$ and $\lambda$ a symmetric affine connection on $U$. The embedding $\iota$ preserves the action by affine combinations on $U_2$ and $V_2$ induced by $\lambda$, respectively, its restriction along $\iota$.  
	\end{lemma}
	\begin{proof}
	    $(i)\quad$~We begin with deriving the familiar transformation law for Christoffel symbols. Let $\la P,Q\ra\in V\la 2\ra$ and $\la P,S\ra\in V\la 2\ra$ be first-order neighbours in $V$. Let $\tilde{\Gamma}_P$ denote the Christoffel symbol of the restriction of the connection $\lambda$ to $V$ along $\iota$ at point $P$. By definition we have
	   $$\iota(Q+S-P+\tilde{\Gamma}_P[Q-P,S-P])=\iota(Q)+\iota(S)-\iota(P)+\Gamma_{\iota(P)}[\iota(Q)-\iota(P),\iota(S)-\iota(P)]$$
	   Due to $\la P,Q,S\ra\in V_2\la 3\ra$ it is
	   $$Q-P+S-P+\tilde{\Gamma}_P[Q-P,S-P]\in D_2(n)$$
	   Taylor-expanding the left hand side yields
	   \begin{multline*}
	       \iota(Q+S-P+\tilde{\Gamma}_P[Q-P,S-P])=\iota(P)+\partial\iota(P)\bigl[Q-P+S-P+\tilde{\Gamma}_P[Q-P,S-P]\bigr] \\
	   +\frac{1}{2}\partial^2\iota(P)\bigl[Q-P+S-P+\tilde{\Gamma}_P[Q-P,S-P]\bigr]^2
	   \end{multline*}
	   Due to $\la P,Q,S\ra\in V_2\la 3\ra$ we find
	   $$ \partial^2\iota(P)\bigl[Q-P+S-P+\tilde{\Gamma}_P[Q-P,S-P]\bigr]^2=\partial^2\iota(P)[Q-P+S-P]^2$$
	   and hence
	   \begin{multline*}
	       \iota(Q+S-P+\tilde{\Gamma}_P[Q-P,S-P])=\iota(P)+\partial\iota(P)\bigl[Q-P+S-P]\\+\frac{1}{2}\partial^2\iota(P)\bigl[Q-P+S-P]^2 +\partial\iota(P)\bigl[\tilde{\Gamma}_P[Q-P,S-P]\bigr]
	   \end{multline*}
	   Since $\la P,Q\ra\in V\la 2\ra$ and $\la P,S\ra\in V\la 2\ra$ further expanding the terms yields
	   \begin{multline*}
	    \iota(P)+\partial\iota(P)\bigl[Q-P+S-P]+\frac{1}{2}\partial^2\iota(P)\bigl[Q-P+S-P]^2\\
	    =\iota(P)+\partial\iota(P)\bigl[Q-P]+\iota(P)+\partial\iota(P)\bigl[S-P]-\iota(P)+\partial^2\iota(P)\bigl[Q-P,S-P]
	   \end{multline*}
	    This simplifies to
	   $$\iota(Q)+\iota(S)-\iota(P)+\partial^2\iota(P)\bigl[Q-P,S-P]$$
	   and finally yields the well-known transformation law of Christoffel symbols
	   $$ \Gamma_{\iota(P)}[\iota(Q)-\iota(P),\iota(S)-\iota(P)]=\partial\iota(P)\bigl[\tilde{\Gamma}_P[Q-P,S-P]\bigr]+\partial^2\iota(P)\bigl[Q-P,S-P]$$
	   $(ii)\quad$~In the second step we apply the same techniques together with this formula to the action of affine combinations defined above. Let $\la P_1,\ldots,P_n\ra\in V_2\la n\ra$
	   $$\iota(\lambda\cdot\la P_1,\ldots,P_n\ra)=\iota\bigl(\sum_{j=1}^n\lambda_jP_j + \frac{1}{2}\bigl(\tilde{\Gamma}_{P_1}\bigl[\sum_{j=1}^n\lambda_jP_j-P_1\bigr]^2-\sum_{j=1}^n\lambda_j\tilde{\Gamma}_{P_1}[P_j-P_1]^2\bigr)\bigr)$$
	   After Taylor-expansion and simplification of the $\partial^2\iota(P)$-term as in step~$(i)$ we get
	   \begin{multline*}
	       \iota(\lambda\cdot\la P_1,\ldots,P_n\ra)=\iota(P_1)+\partial\iota(P_1)\bigl[\sum_{j=1}^n\lambda_j P_j-P_1\bigr] +\frac{1}{2}\partial^2\iota(P_1)\bigl[\sum_{j=1}^n\lambda_j P_j-P_1\bigr]^2\\ +\frac{1}{2}\bigl(\partial\iota(P_1)\bigl[\tilde{\Gamma}_{P_1}\bigl[\sum_{j=1}^n\lambda_jP_j-P_1\bigr]^2\bigr]-\sum_{j=1}^n\lambda_j\partial\iota(P_1)\bigl[\tilde{\Gamma}_{P_1}[P_j-P_1]^2\bigr]\bigr)
	   \end{multline*}
	   Applying the transformation law of the Christoffel symbols yields
	   \begin{align*}
	       \iota(\lambda\cdot\la P_1,\ldots,P_n\ra) =&\iota(P_1)+\partial\iota(P_1)\bigl[\sum_{j=1}^n\lambda_j P_j-P_1\bigr] 
	       +\sum_{j=1}^n\lambda_j\frac{1}{2}\partial^2\iota(P_1)\bigl[P_j-P_1\bigr]^2 \\
	       &+\frac{1}{2}\bigl(\Gamma_{\iota(P_1)}\bigl[\sum_{j=1}^n\lambda_j\iota(P_j)-\iota(P_1)\bigr]^2-\sum_{j=1}^n\lambda_j\Gamma_{\iota(P_1)}[\iota(P_j)-\iota(P_1)]^2\bigr)
	   \end{align*}
	   Using $\sum_{j=1}^n\lambda_j=1$ we find 
	   \begin{align*}
	       \iota(P_1)+\partial\iota(P_1)\bigl[\sum_{j=1}^n\lambda_j P_j-P_1\bigr] 
	       &+ \sum_{j=1}^n\lambda_j\frac{1}{2}\partial^2\iota(P_1)[P_j-P_1]^2 \\
	       &= \iota(P_1)+\sum_{j=1}^n\lambda_j\bigl(\partial\iota(P_1)[ P_j-P_1]
	       +\frac{1}{2}\partial^2\iota(P_1)[P_j-P_1]^2\bigr) \\
	       &=\iota(P_1)+\sum_{j=1}^n\lambda_j(\iota(P_j)-\iota(P_1)) \\
	       &=\sum_{j=1}^n\lambda_j\iota(P_j)
	   \end{align*}
	    Substituting this in the equation above yields the desired
	    \begin{align*}
	        \iota(\lambda\cdot\la P_1,\ldots,P_n\ra) &=\sum_{j=1}^n\lambda_j \iota(P_j) +\frac{1}{2}\bigl(\Gamma_{\iota(P_1)}\bigl[\sum_{j=1}^n\lambda_j\iota(P_j)-\iota(P_1)\bigr]^2-\sum_{j=1}^n\lambda_j\Gamma_{\iota(P_1)}[\iota(P_j)-\iota(P_1)]^2\bigr) \\
	        &= \lambda\cdot\la\iota(P_1),\ldots,\iota(P_n)\ra
	    \end{align*}
	\end{proof}
    It remains to show that the action of affine combinations on $U_2$ satisfies the associativity axiom and hence is a 2nd-order i-affine structure. This follows from another lengthy calculation following the same techniques we have been using above: Taylor-expansion and vanishing of terms which are $k$-linear for $k\geq 3$. We shall only give the main steps.
    \begin{lemma} Let $\lambda$ be a symmetric affine connection on a formally open subspace $U$ of some $R^N$. Let $\Gamma$ denote the Christoffel symbol of $\lambda$. The action of $\calA(n)$ on $U_2\la n\ra$ defined by
	$$\mu\cdot\la P_1,\ldots,P_n\ra = \sum_{j=1}^n\mu_jP_j + \frac{1}{2}\bigl(\Gamma_{P_1}\bigl[\sum_{j=1}^n\mu_jP_j-P_1\bigr]^2-\sum_{j=1}^n\mu_j
	\Gamma_{P_1}[P_j-P_1]^2\bigr)$$
    for each $n\geq 1$ defines a 2nd-order i-affine structure on $U$.
    \end{lemma}
    \begin{proof}
        It remains to show the associativity axiom, i.e. for all $\lambda^1,\ldots,\lambda^m\in\calA(n)$ and $\mu\in\calA(m)$ we have
        $$\mu\cdot\bigl\la\lambda^1\cdot\la P_1,\ldots,P_n\ra,\ldots,\lambda^m\cdot\la P_1,\ldots,P_n\ra\bigr\ra =
        \bigl(\sum_{\ell=1}^m\mu_\ell\lambda^\ell\bigr)\cdot\la P_1,\ldots,P_n\ra$$
        The right hand side is by definition
        $$\sum_{j=1}^n\bigl(\sum_{\ell=1}^m\mu_\ell\lambda^\ell_j\bigr)P_j + \frac{1}{2}\bigl(\Gamma_{P_1}\bigl[\sum_{j=1}^n\bigl(\sum_{\ell=1}^m\mu_\ell\lambda^\ell_j\bigr)P_j-P_1\bigr]^2-\sum_{j=1}^n\bigl(\sum_{\ell=1}^m\mu_\ell\lambda^\ell_j\bigr)\Gamma_{P_1}[P_j-P_1]^2\bigr)$$
        The left hand side evaluates to
        \begin{multline*}
            \sum_{\ell=1}^m\mu_\ell\;\lambda^\ell\cdot\la P_1,\ldots,P_n\ra + \frac{1}{2}\bigl(\Gamma_{\lambda^1\cdot\la P_1,\ldots,P_n\ra}\bigl[\sum_{\ell=1}^m\mu_\ell\;\lambda^\ell\cdot\la P_1,\ldots,P_n\ra-\lambda^1\cdot\la P_1,\ldots,P_n\ra\bigr]^2\\ -\sum_{\ell=1}^m\mu_\ell\Gamma_{\lambda^1\cdot\la P_1,\ldots,P_n\ra}[\lambda^\ell\cdot\la P_1,\ldots,P_n\ra-\lambda^1\cdot\la P_1,\ldots,P_n\ra]^2\bigr)
        \end{multline*}
        Evaluating the first term yields
        \begin{multline*}
            \sum_{\ell=1}^m\mu_\ell\;\lambda^\ell\cdot\la P_1,\ldots,P_n\ra=
            \sum_{j=1}^n\bigl(\sum_{\ell=1}^m\mu_\ell\lambda^\ell_j\bigr)P_j + \frac{1}{2}\bigl(\sum_{\ell=1}^m\mu_\ell\Gamma_{P_1}\bigl[\sum_{j=1}^n\lambda^\ell_jP_j-P_1\bigr]^2\\-\sum_{j=1}^n\bigl(\sum_{\ell=1}^m\mu_\ell\lambda^\ell_j\bigr)\Gamma_{P_1}[P_j-P_1]^2\bigr)
        \end{multline*}
        Comparing this with the right hand side of the associativity condition reveals that for the latter to hold we need to show 
        \begin{multline*}
            \Gamma_{P_1}\bigl[\sum_{j=1}^n\bigl(\sum_{\ell=1}^m\mu_\ell\lambda^\ell_j\bigr)P_j-P_1\bigr]^2 \\= \Gamma_{\lambda^1\cdot\la P_1,\ldots,P_n\ra}\bigl[\sum_{\ell=1}^m\mu_\ell\;\lambda^\ell\cdot\la P_1,\ldots,P_n\ra-\lambda^1\cdot\la P_1,\ldots,P_n\ra\bigr]^2 +\sum_{\ell=1}^m\mu_\ell\bigl(\Gamma_{P_1}\bigl[\sum_{j=1}^n\lambda^\ell_jP_j-P_1\bigr]^2\\-\Gamma_{\lambda^1\cdot\la P_1,\ldots,P_n\ra}[\lambda^\ell\cdot\la P_1,\ldots,P_n\ra-\lambda^1\cdot\la P_1,\ldots,P_n\ra]^2\bigr)
        \end{multline*}
        Due to $\la P_1,\ldots,P_n\ra\in U_2\la n\ra$ the Christoffel symbols simplify to
        \begin{align*}
            &\Gamma_{\lambda^1\cdot\la P_1,\ldots,P_n\ra}\bigl[\sum_{\ell=1}^m\mu_\ell\;\lambda^\ell\cdot\la P_1,\ldots,P_n\ra-\lambda^1\cdot\la P_1,\ldots,P_n\ra\bigr]^2 \\
            &\hskip 18ex=\Gamma_{\lambda^1\cdot\la P_1,\ldots,P_n\ra}\bigl[\sum_{j=1}^n\bigl(\sum_{\ell=1}^m\mu_\ell\lambda^\ell_j\bigr)P_j-\sum_{j=1}^n\lambda^1_j P_j\bigr]^2\\
            &\Gamma_{\lambda^1\cdot\la P_1,\ldots,P_n\ra}[\lambda^\ell\cdot\la P_1,\ldots,P_n\ra-\lambda^1\cdot\la P_1,\ldots,P_n\ra]^2\\
            &\hskip 18ex =\Gamma_{\lambda^1\cdot\la P_1,\ldots,P_n\ra}\bigl[\sum_{j=1}^n\lambda^\ell_j P_j-\sum_{j=1}^n\lambda^1_j P_j\bigr]^2
        \end{align*}
        Furthermore, $\lambda^1\cdot\la P_1,\ldots,P_n\ra-P_1\in D_2(n)$ and the Taylor-expansion of $\Gamma_{\lambda^1\cdot\la P_1,\ldots,P_n\ra}$ at $P_1$ yields
        \begin{align*}
            \Gamma_{\lambda^1\cdot\la P_1,\ldots,P_n\ra}\bigl[\sum_{j=1}^n\bigl(\sum_{\ell=1}^m\mu_\ell\lambda^\ell_j\bigr)P_j-\sum_{j=1}^n\lambda^1_j P_j\bigr]^2 &= \Gamma_{P_1}\bigl[\sum_{j=1}^n\bigl(\sum_{\ell=1}^m\mu_\ell\lambda^\ell_j\bigr)P_j-\sum_{j=1}^n\lambda^1_j P_j\bigr]^2 \\
            \Gamma_{\lambda^1\cdot\la P_1,\ldots,P_n\ra}\bigl[\sum_{j=1}^n\lambda^\ell_j P_j-\sum_{j=1}^n\lambda^1_j P_j\bigr]^2 &= \Gamma_{P_1}\bigl[\sum_{j=1}^n\lambda^\ell_j P_j-\sum_{j=1}^n\lambda^1_j P_j\bigr]^2
        \end{align*}
        since all the other terms contain $k$-linear occurrences of $P_j-P_i$ for $k\geq 3$ and thus vanish. Therefore, it is sufficient to show the equation
        \begin{multline*}
            \Gamma_{P_1}\bigl[\sum_{j=1}^n\bigl(\sum_{\ell=1}^m\mu_\ell\lambda^\ell_j\bigr)P_j-P_1\bigr]^2 = \Gamma_{P_1}\bigl[\sum_{j=1}^n\bigl(\sum_{\ell=1}^m\mu_\ell\lambda^\ell_j\bigr)P_j-\sum_{j=1}^n\lambda^1_j P_j\bigr]^2\\ +\sum_{\ell=1}^m\mu_\ell\bigl(\Gamma_{P_1}\bigl[\sum_{j=1}^n\lambda^\ell_jP_j-P_1\bigr]^2-\Gamma_{P_1}\bigl[\sum_{j=1}^n\lambda^\ell_j P_j-\sum_{j=1}^n\lambda^1_j P_j\bigr]^2\bigr)
        \end{multline*}
        By adding $-P_1+P_1$ to the first argument of the second $\Gamma_{P_1}$ and using symmetric bilinearity we find 
        \begin{align*}
           \Gamma_{P_1}\bigl[\sum_{j=1}^n\lambda^\ell_jP_j&-P_1\bigr]^2-\Gamma_{P_1}\bigl[\sum_{j=1}^n\lambda^\ell_j P_j-\sum_{j=1}^n\lambda^1_j P_j\bigr]^2\\
           &=\Gamma_{P_1}\bigl[\sum_{j=1}^n\lambda^\ell_jP_j-P_1\bigr]^2 - \Gamma_{P_1}\bigl[\sum_{j=1}^n\lambda^\ell_j P_j-P_1,\sum_{j=1}^n\lambda^\ell_j P_j-\sum_{j=1}^n\lambda^1_j P_j\bigr]\\
           &\quad + \Gamma_{P_1}\bigl[\sum_{j=1}^n\lambda^1_j P_j-P_1,\sum_{j=1}^n\lambda^\ell_j P_j-\sum_{j=1}^n\lambda^1_j P_j\bigr]\\
           &= \Gamma_{P_1}\bigl[2\sum_{j=1}^n\lambda^\ell_jP_j-P_1-\sum_{j=1}^n\lambda^1_j P_j,\sum_{j=1}^n\lambda^1_j P_j-P_1\bigr]\\
            &=2\Gamma_{P_1}\bigl[\sum_{j=1}^n\lambda^\ell_jP_j-P_1,\sum_{j=1}^n\lambda^1_j P_j-P_1\bigr] -\Gamma_{P_1}\bigl[\sum_{j=1}^n\lambda^1_j P_j-P_1\bigr]^2
        \end{align*}
        and hence
        \begin{multline*}
            \sum_{\ell=1}^m\mu_\ell\bigl(\Gamma_{P_1}\bigl[\sum_{j=1}^n\lambda^\ell_jP_j-P_1\bigr]^2-\Gamma_{P_1}\bigl[\sum_{j=1}^n\lambda^\ell_j P_j-\sum_{j=1}^n\lambda^1_j P_j\bigr]^2\bigr) \\  =2\Gamma_{P_1}\bigl[\sum_{j=1}^n\bigl(\sum_{\ell=1}^m\mu_\ell\lambda^\ell_j\bigr)P_j-P_1,\sum_{j=1}^n\lambda^1_j P_j-P_1\bigr]
               -\Gamma_{P_1}\bigl[\sum_{j=1}^n\lambda^1_j P_j-P_1\bigr]^2
        \end{multline*}
        As regards the first term on the right hand side of the equation we wish to show, we find
        \begin{multline*}
           \Gamma_{P_1}\bigl[\sum_{j=1}^n\bigl(\sum_{\ell=1}^m\mu_\ell\lambda^\ell_j\bigr)P_j-\sum_{j=1}^n\lambda^1_j P_j\bigr]^2 = \Gamma_{P_1}\bigl[\sum_{j=1}^n\bigl(\sum_{\ell=1}^m\mu_\ell\lambda^\ell_j\bigr)P_j-P_1\bigr]^2 + \Gamma_{P_1}\bigl[\sum_{j=1}^n\lambda^1_j P_j-P_1\bigr]^2 \\- 2\Gamma_{P_1}\bigl[\sum_{j=1}^n\bigl(\sum_{\ell=1}^m\mu_\ell\lambda^\ell_j\bigr)P_j-P_1,\sum_{j=1}^n\lambda^1_j P_j-P_1\bigr]
        \end{multline*}
        and thus
        \begin{multline*}
             \Gamma_{P_1}\bigl[\sum_{j=1}^n\bigl(\sum_{\ell=1}^m\mu_\ell\lambda^\ell_j\bigr)P_j-\sum_{j=1}^n\lambda^1_j P_j\bigr]^2 \\ +\sum_{\ell=1}^m\mu_\ell\bigl(\Gamma_{P_1}\bigl[\sum_{j=1}^n\lambda^\ell_jP_j-P_1\bigr]^2-\Gamma_{P_1}[\sum_{j=1}^n\lambda^\ell_j P_j-\sum_{j=1}^n\lambda^1_j P_j]^2\bigr)
             =\Gamma_{P_1}\bigl[\sum_{j=1}^n\bigl(\sum_{\ell=1}^m\mu_\ell\lambda^\ell_j\bigr)P_j-P_1\bigr]^2
        \end{multline*}
        as required.
    \end{proof}
    
	\begin{theorem}\label{thm:equiv_iaff_conn} 
		Every symmetric affine connection $\lambda$ on a formal manifold $A$ extends to an i-affine structure on $A_2$ in such a way that 
		$$
			\lambda(P,Q,S) = (-1,1,1)\cdot\la P,Q,S\ra
		$$
		for all $(P,Q,S)\in A^3$ such that $\la P,Q\ra, \la P,S\ra\in A_1\la 2\ra$. (The right hand side denotes the i-affine combination induced by $\lambda$ on $A_2$ as defined above.) 
	\end{theorem}
	\begin{proof}
		It remains to show that $\lambda$ agrees with the given affine combination of the induced 2nd-order i-structure. As shown in proposition~\ref{prop:2nd-iaff-connection} we have $\la P,Q,S\ra\in A_2\la 3\ra$. We consider everything in a chart $U$. By definition we have
		$$
			(-1,1,1)\cdot\la P,Q,S\ra = -P + Q + S +\frac{1}{2}(\Gamma_P[Q-P + S-P]^2 - \Gamma_P[Q-P]^2 - \Gamma_P[S-P]^2)
		$$
		Expanding the symmetric bilinear map $\Gamma_P$ results in
		$$
			(-1,1,1)\cdot\la P,Q,S\ra = \lambda(P,Q,S)
		$$
		as claimed.
		%
%
	\end{proof}
	\section{Existence results for 2nd-order i-affine structures}
	
	It remains to show that a manifold admits a 2nd-order i-affine structure. Due to theorem~\ref{thm:equiv_iaff_conn} this is equivalent to showing that it admits a symmetric affine connection. The author is not aware of an existence result of affine connections of points on a \emph{formal} manifold. However, for a \emph{smooth} manifold $A$ (considered as being embedded in a well-adapted model of SDG) there are various ways to show the existence of an affine connection on points. For example, one can use that every smooth manifold admits a Riemannian metric and construct a Levi-Civita connection on points \cite{Riem}. Combining this with theorem~\ref{thm:equiv_iaff_conn} yields the first existence result:
	
	\begin{corollary}\label{cor:exist_2nd_i-aff} Every smooth manifold admits a 2nd-order i-affine structure.
	\end{corollary}
	
	\begin{remark}\label{rem:metric} Even though a Riemannian metric is classically defined as a positive definite symmetric bilinear form on the tangent vectors of a manifold, it is possible to also construct Riemannian metrics on points as defined in \cite{Riem}, \cite[chapter~8]{SGM}. For smooth manifolds in a \emph{well-adapted model} the most direct way is to proceed as in the classical construction of a Riemannian metric on the tangent bundle: We use the existence of a locally finite atlas and its subordinated partition of unity to glue together the canonical metrics on the formally open subsets $U\into R^n$
		$$
			g: U_2\la 2\ra \to R,\qquad (P,Q)\mapsto (Q-P)\bullet(Q-P)
		$$  
	induced by the scalar product of $R^n$.	Another construction, which seems to be more in the spirit of SDG, is to use the $\log$-$\exp$-bijection as discussed in chapter~4.3 in \cite{SGM} to show that a Riemannian metric on tangent vectors induces a Riemannian metric on points. This, however, only defines a metric on 2nd-order neighbours $(P,Q)$ with the property that there is a point $X$, such that $(P,X)$ and $(X,Q)$ are 1st-order neighbours:
	$$
		g(P,Q)=g_X(\log_X(P),\log_X(Q))
	$$
	where the metric on the right hand side denotes the Riemannian metric on the tangent space at $X$. (See \cite[chapters 4.3, 8]{SGM} for definitions.) Although such a $g$ has a unique extension to $U_2\la 2\ra$ in a chart $U$, to obtain an extension to the whole 2nd neighbourhood of the diagonal of the smooth manifold seems to require the existence of an $O(n)$-atlas.
	   
	\end{remark}
	
	For a smooth (regular) \emph{submanifold} $M$ of $R^n$ there is another way to construct a second-order i-affine structure inspired by the following construction: a symmetric affine connection on $M$ can be obtained by applying the canonical (flat) connection of $R^n$ to two tangent vectors over the same base point $P$ and then project the resulting vector back into the tangent space $T_PM$ along its normal space $N_PM$. 
	
	Since we are interested in points instead of tangent vectors we ought to replace the projection on the vector bundle with a mapping on the base spaces. This is possible due to the \emph{tubular neighbourhood theorem}, from which one can conclude that every submanifold $M$ is a retract of a (formally) open subspace $U\subset R^n$ (see theorem~5.1 in \cite{Hirsch}, for example). (Here we are making use of well-adapted models once more, and the fact that the embedding maps open subsets of $\IR^n$ to formally open subspaces of $R^n$; see \cite{WAM} but also theorem~III.3.4 in \cite{SDG}.)
	
	We shall thus consider a retract $\iota: M\into U\subseteq R^n$ of a formally open subspace $U$ with retraction $r: U\to M$. Pulling back the 2nd-order i-structure $U_2$ via $\iota$ yields an i-structure on $M$, which we shall denote by $M_2$. We wish to define an i-affine structure on $M_2$ by projecting the 2nd-order i-affine structure on $U_2$ via $r$. For each $n\geq 1$ we define an action of $\calA(n)$ on $M_2\la n\ra$ by
	$$
		\lambda\cdot\la P_1,\ldots,P_n\ra = r\Bigl(\sum_{j=1}^n \lambda_j \iota(P_j)\Bigr)
	$$
	By theorem~\ref{thm:kth-i-morph} the \emph{idempotent} $e=\iota \circ r$ is an i-morphism. Due to the construction of $M_2$ we can conclude that $r$ is an i-morphism. From this and the properties of the i-affine structure on $U_2$ it follows readily that the action defined above satisfies the neighbourhood and projection axioms. It remains to show the associativity axiom.
	
	Let $\lambda^1,\ldots,\lambda^m\in\calA(n)$, $\mu\in\calA(m)$ and $\la P_1,\ldots,P_n\ra\in M_2\la n\ra$ for some $m,n\geq 1$. We need to show that
	$$
		r\Bigl(\sum_{\ell=1}^m \mu_\ell\, e\bigl(\sum_{j=1}^n\lambda^\ell_j \iota(P_j)\bigr)\Bigr) = r\Bigl(\sum_{\ell=1}^m\sum_{j=1}^n \mu_\ell\lambda^\ell_j \iota(P_j)\Bigr) 
	$$
	To lighten the notation for the subsequent calculations we shall identify $A$ with its image $\iota(A)$. In this case the equivalent equation we have to show is
	$$
		e\Bigl(\sum_{\ell=1}^m \mu_\ell\, e\bigl(\sum_{j=1}^n\lambda^\ell_j P_j\bigr)\Bigr) = e\Bigl(\sum_{\ell=1}^m\sum_{j=1}^n \mu_\ell\lambda^\ell_j P_j\Bigr) 
	$$  
	The rest of the proof is once again a direct calculation based on Taylor-expansion combined with multilinear algebra of nil-potents. 
	
	Due to $\sum_{j=1}^n \lambda^\ell_j=1$ we can write 
	$$
		\sum_{j=1}^n \lambda^\ell_j P_j = P_1 + \sum_{j=1}^n \lambda^\ell_j(P_j-P_1)
	$$
	Since $\la P_1,\ldots,P_n\ra\in M_2\la n\ra$ we have $\sum_{j=1}^n \lambda^\ell_j(P_j-P_1)\in D_2(n)$. The first Taylor expansion yields (note that $e(P_j)=P_j$)
	%
	%
	$$
		e\Bigl(\sum_{j=1}^n \lambda^\ell_j P_j\Bigr) = P_1 + \partial e(P_1)\Bigl[\sum_{j=1}^n \lambda^\ell_j (P_j-P_1)\Bigr] + \frac{1}{2}\partial^2 e(P_1)\Bigl[\sum_{j=1}^n \lambda^\ell_j (P_j-P_1)\Bigr]^2
	$$ 
	and hence
	$$
		\sum_{\ell=1}^{m}\mu_\ell\, e\Bigl(\sum_{j=1}^n \lambda^\ell_j P_j\Bigr) = P_1 + \partial e(P_1)\Bigl[\sum_{\ell=1}^{m}\sum_{j=1}^n \mu_\ell\lambda^\ell_j (P_j-P_1)\Bigr] + \sum_{\ell=1}^{m}\mu_\ell\frac{1}{2}\partial^2 e(P_1)\Bigl[\sum_{j=1}^n \lambda^\ell_j (P_j-P_1)\Bigr]^2
	$$ 
	Omitting all the vanishing $k$-linear terms in $P_j-P_1$ for $k\geq 3$ the Taylor expansion of $e(\sum_{\ell=1}^{m}\mu_\ell\, e(\sum_{j=1}^n \lambda^\ell_j P_j))$ at $P_1$ reads
	\begin{align*}
		e\Bigl(\sum_{\ell=1}^{m}\mu_\ell\, e\bigl(\sum_{j=1}^n \lambda^\ell_j P_j\bigr)\Bigr) =  P_1 & + (\partial e(P_1))^2\Bigl[\sum_{\ell=1}^{m}\sum_{j=1}^n\mu_\ell \lambda^\ell_j (P_j-P_1)\Bigr] \\
		& + \sum_{\ell=1}^{m}\mu_\ell\frac{1}{2}\partial e(P_1)\circ\partial^2 e(P_1)\Bigl[\sum_{j=1}^n \lambda^\ell_j (P_j-P_1)\Bigr]^2 \\
		& +  \frac{1}{2}\partial^2 e(P_1)\Bigl[\partial e(P_1)\bigl[\sum_{\ell=1}^{m}\sum_{j=1}^n\mu_\ell \lambda^\ell_j (P_j-P_1)\bigr]\Bigr]^2
	\end{align*}
	To show that the associativity axiom holds we shall show that the right hand side simplifies to the Taylor expansion of $e(\sum_{\ell=1}^{m}\sum_{j=1}^n \mu_\ell\lambda^\ell_j P_j)$ at $P_1$:
	$$
		e\Bigl(\sum_{\ell=1}^{m}\sum_{j=1}^n \mu_\ell\lambda^\ell_j P_j\Bigr) = P_1 + \partial e(P_1)\Bigl[\sum_{\ell=1}^{m}\sum_{j=1}^n \mu_\ell\lambda^\ell_j (P_j-P_1)\Bigr] + \frac{1}{2}\partial^2 e(P_1)\Bigl[\sum_{\ell=1}^{m}\sum_{j=1}^n \mu_\ell\lambda^\ell_j (P_j-P_1)\Bigr]^2
	$$
	Differentiating $e^2=e$ at $P\in U$ yields
	%
	$$
		\partial e(e(P))\circ \partial e(P) = \partial e(P)
	$$
	Differentiating a second time yields
	$$
		\partial^2 e(e(P))[\partial e(p)]^2 + \partial e(e(P))\circ \partial^2 e(P) = \partial^2 e(P)
	$$
	If $P\in M$, i.e. $e(P)=P$, this simplifies to $(\partial e(P))^2=\partial e(P)$ and
	$$
		\partial^2 e(P)[\partial e(P)]^2 + \partial e(P)\circ \partial^2 e(P) = \partial^2 e(P)
	$$
	Since
	$$
		P_j - P_1 = e(P_j) - e(P_1) = \partial e(P_1)[P_j-P_1] + \frac{1}{2}\partial^2 e(P_1)[P_j-P_1]^2
	$$
	and $\la P_1, P_j\ra \in U_2\la 2\ra$ we have
	$$
		\partial^2 e(P_1)\bigl[\partial e(P_1)[P_j-P_1]\bigr]^2 = \partial^2 e(P_1)[P_j-P_1]^2
	$$
	Substituting this into the equation obtained from differentiating $e^2=e$ twice yields
	$$
		\partial e(P_1)\circ \partial^2 e(P_1)[P_j-P_1]^2 = 0 
	$$
	and thus
	$$
		\sum_{\ell=1}^{m}\mu_\ell\frac{1}{2}\partial e(P_1)\circ\partial^2 e(P_1)\Bigl[\sum_{j=1}^n \lambda^\ell_j (P_j-P_1)\Bigr]^2 = 0
	$$
	Moreover, we find
	$$
		\frac{1}{2}\partial^2 e(P_1)\Bigl[\partial e(P_1)\bigl[\sum_{\ell=1}^{m}\sum_{j=1}^n\mu_\ell \lambda^\ell_j (P_j-P_1)\bigr]\Bigr]^2 = \frac{1}{2}\partial^2 e(P_1)\Bigl[\sum_{\ell=1}^{m}\sum_{j=1}^n\mu_\ell \lambda^\ell_j (P_j-P_1)\Bigr]^2
	$$
	Finally, applying $(\partial e(P_1))^2 = \partial e(P_1)$ results in
	$$
		e\Bigl(\sum_{\ell=1}^{m}\mu_\ell\, e\bigl(\sum_{j=1}^n \lambda^\ell_j P_j\bigr)\Bigr) =  P_1  + \partial e(P_1)\Bigl[\sum_{\ell=1}^{m}\sum_{j=1}^n\mu_\ell \lambda^\ell_j (P_j-P_1)\Bigr]
		 +  \frac{1}{2}\partial^2 e(P_1)\Bigl[\sum_{\ell=1}^{m}\sum_{j=1}^n\mu_\ell \lambda^\ell_j (P_j-P_1)\Bigr]^2	
	$$
	which is equal to $e(\sum_{\ell=1}^{m}\sum_{j=1}^n \mu_\ell\lambda^\ell_j P_j)$ as claimed. This concludes showing that the action we have defined above is indeed an i-affine structure on $M_2$.
	
	\begin{theorem}\label{thm:2nd-order-iaff-on-retracts}
		Let $U\subseteq R^n$ be formally open and $\iota: M\into U$ a retract with retraction $r: U\to M$. Defining $M_2$ as the pullback of the 2nd-order i-structure $U_2$ along $\iota$ and setting 
		$$
		\lambda\cdot\la P_1,\ldots,P_n\ra = r\Bigl(\sum_{j=1}^n \lambda_j \iota(P_j)\Bigr)
		$$
		for each $n\geq 1$, $\lambda\in \calA(n)$ and $\la P_1,\ldots,P_n\ra\in M_2\la n\ra$ makes $M$ into an i-affine space.
	\end{theorem}
	
	\begin{remark}
		\begin{enumerate}[(a)]
			\item In the case that $M$ is a smooth (regular) submanifold we note that the 2nd-order i-structure defined via the charts in theorem~\ref{thm:kth-i-struc} agrees with the pullback of the 2nd-order i-structure via $\iota$. Indeed, also for the 2nd-order i-structure in theorem~\ref{thm:kth-i-struc} the map $\iota$ is an i-structure \emph{reflecting} i-morphism\footnote{A map $f:X\to Y$ is said to \emph{reflect i-structure}, if $\la f(P_1),\ldots, f(P_n)\ra\in Y\la n\ra$ implies $\la P_1,\ldots, P_n\ra\in X\la n\ra$ for all $n\geq 1$.}, so both i-structures agree. 
			This follows from considering submanifold charts of $M$, i.e. pullbacks (of certain) $U$-charts $V$ along $\iota$. By construction $U$-charts $V$ reflect the 2nd-order i-structure of $R^n$ and submanifold charts $V\cap M$ reflect the i-structure of $M$. The composite is an embedding of $V\cap M$ into a subspace $R^m\times\{0\}$ of $R^n$, which is readily seen to reflect the 2nd-order i-structure as $V\cap M$ is formally open in $R^m$ and the embedding $R^m\cong R^m\times\{0\}\into R^n$ reflects i-structure.  
			%
			\item For a formal manifold $M$ the proof can be simplified. Due to theorem~\ref{thm:equiv_iaff_conn} it suffices to notice that $\lambda(P,Q,S) = r(\iota(Q) + \iota(S) - \iota(P))$ is a symmetric affine connection on $M$. Although the proof presented here works for more general spaces than manifolds, we are unable to provide a non-algebraic example of such. 
			(Despite there being more general i-affine spaces than smooth manifolds like Euclidean neighbourhood retracts, for example, the morphisms between these spaces are smooth maps, which means that the property of being a retract only holds in a  well-adapted model when the retraction is smooth; but this is the case if and only if the retract is a submanifold.)  
			%
			\item We could have also considered 2nd-order i-affine structures on $U_2$ different from the canonical one induced by $R^n$ to construct the 2nd-order i-affine structure on $M_2$. In the case of $M$ being a manifold this is obvious as we only need a symmetric affine connection as pointed out in the previous remark. However, by making use of Christoffel symbols the proof by direct calculation should generalise to this case as well (although it will be far more tedious).  
		\end{enumerate}
	\end{remark}    
	
    \section{Conclusion}\label{sec:conclusion}
	
	An action of (the clone of) affine combinations on an i-structure is an algebraic model that makes precise the long-standing idea of differential geometry and of calculus that a (smooth) space has a geometry that is affine at the infinitesimal scale. These algebraic structures have been extracted by the author from Kock's work \cite{SDG}, \cite{SGM}. The author has then generalised and studied them as infinitesimal models of algebraic theories in \cite{IMAT}. 
	
	Within the framework of Synthetic Differential Geometry, in particular within the algebraic and well-adapted models of SDG there is a wealth of examples of i-affine spaces besides that of smooth and formal manifolds. This means that \emph{the same} infinitesimal constructs and \emph{the same} algebra of infinitesimals can be applied much more widely and beyond the context of (smooth) manifolds. However, so far (almost) all these examples have been based on the nil-square i-structure only\footnote{The only exception has been the pointwise i-affine structure on function spaces studied in \cite[chap.~3.3.2]{IMAT}.}. 
	
	In this paper we have shown that besides the canonical nil-square i-structure, a formal manifold carries a natural $k$th-order i-structure for each $k\geq 1$. The affine structure on $R^n$ induces i-affine structures on each of its $k$th-order i-structures. In contrast to the nil-square i-affine structure the $k$th-order i-affine structures for $k\geq 2$ are not preserved by all maps $R^n\to R^m$, and are hence not natural anymore. However, as we have shown for formal manifolds, there is a correspondence between symmetric affine connections (on points) and 2nd-order i-affine structures. This provides us with a first example that a higher-order i-affine structure can be obtained from the data of a higher-order geometric structure on a formal manifold.    
	     
	It should be noticed that except for the gluing arguments all the proofs in this paper involved only Taylor expansions and the multilinear algebra of infinitesimals. This agrees with the common practise when considering (infinitely) small perturbations in physics and engineering; the only difference being that higher-order terms are neglected, but do not vanish. Often the decision which terms are to be neglected is based on intuition and experience rather than on prescribed algebraic rules. However, the presence of such rules allows for a more rigorous and systematic analysis. Although the calculations involving higher-order infinitesimals tend to become lengthy and tedious, they remain purely algebraic manipulations that can be implemented as a module of a computer algebra system. 
	
	\paragraph{Does a manifold admit higher-order i-affine structures?} By showing that each symmetric affine connection induces a second-order i-affine structure on a formal manifold we have a positive answer in the case $k=2$. 
	As regards $k\geq 3$ we can attempt to use the construction in theorem~\ref{thm:2nd-order-iaff-on-retracts} and adapt the calculations. The author was able to carry this out successfully for the case $k=3$; theorem~\ref{thm:2nd-order-iaff-on-retracts} can therefore be generalised to 3rd-order i-affine structures. It remains to be seen whether a generalisation to higher (and possibly all) $k$ is possible as well.
	 
	A related problem of the extension of flat symmetric (pointwise) affine connections to an affine structure has been studied by Kock in \cite{Int_of_conn}. Although Kock has not been using i-affine structures explicitly, his \cite[theorem~3.9]{Int_of_conn} can be reformulated with such: \textit{Every flat symmetric (pointwise) affine connection on a formal manifold $A$ induces an i-affine structure on the $\infty$-i-structure of $A$.} The $\infty$-i-structure is the join (= colimit) of all $k$th-order i-structures for $k\geq 1$. As each $k$th order i-structure is \emph{natural} on a formal manifold, so is the $\infty$-i-structure. However, in the context of \cite{Int_of_conn} we have to restrict the spaces $D_k(n)$ to the set of elements of $R^n$ that can be written as a sum of $k$ first-order infinitesimals $d_j\in D(n)$. Note that this re-definition does not affect any of the proofs given in this paper.    
		
	\paragraph{Can we extend the $k$th-order i-structures without compromising their naturality?} Even though the extension of a connection to a 2nd-order i-affine structure is conceptually satisfying, it might not be too useful in practice: it is not easy to show that a family of points constitutes a 2nd-order i-tuple and the 2nd-order i-structure does not contain the nil-square i-structure that is much easier to work with in this respect. Moreover, one typically arrives at higher-order structures by concatenating successive first-order steps. \emph{Is there a class of multilinear forms for which a 2nd-order i-structure would have one or both of these desirable properties?} The proof of theorem~\ref{thm:kth-i-aff} would work for any class of $(k+1)$-linear forms, but the other results have to be treated with care.	

    \paragraph{Does a symmetric affine connection determine an i-affine structure uniquely?} We have shown that a symmetric affine connection extends to a 2nd-order i-affine structure on formal manifolds. What we have not addressed is the question whether the 2nd-order i-affine structure is uniquely determined by the connection, or, if not, what structure parametrises the possible freedom of choice.
    
    The author was able to show that in a \emph{well-adapted model} each smooth manifold $A$ carries \emph{only one} i-affine structure on the first-order i-structure $A_1$. Studying the uniqueness of 2nd-order i-affine structures is current work in progress.

\end{document}